\documentclass{article}
\usepackage{amsmath, amssymb, latexsym, amsthm, amsfonts, mathrsfs}

\newtheorem{thm}{Theorem}
\newtheorem{prop}[thm]{Proposition}
\newtheorem{lem}[thm]{Lemma}
\newtheorem{cor}[thm]{Corollary}

\begin{document}

\title{On sequences of positive integers containing no $p$ terms in arithmetic progression}

\author{Goutam Pal\thanks{E-mail: \texttt{gkmstr@gmail.com}}}

\date{11 Nov 2006}
\maketitle

\centerline{Dept. of Mathematics,}
\centerline{RCC-Institute of Information Technology,}
\centerline{South Canal Road, Beliaghata,}
\centerline{Kolkata-700015}

\begin{abstract}
We use topological ideas to show that, assuming the conjecture of Erd\"{o}s \cite{UL} on subsets of positive integers having no $p$ terms in arithmetic progression (A. P.), there must exist a subset $M_p$ of positive integers with no $p$ terms in A. P. with the property that among all such subsets, $M_p$ maximizes the sum of the reciprocals of its elements. 
\end{abstract}

\section{Introduction}

A famous conjecture of Erd\"{o}s  asserts that if $A$ is a
subset of the positive integers having the property that $\sum_{a\in A} \frac{1}{a} =
\infty$, then $A$ must contain arithmetic progressions of arbitrary length. 
A special case of the conjecture, when $A$ is the set of prime numbers, was recently proved by Green and Tao \cite{GT}. This implies that if a subset $A$ of the set of positive integers contains no
arithmetic progression of length $p$, where $p \geq 3$ is a fixed integer, then
the sum $\sum_{a\in A} \frac{1}{a}$ must converge. In this paper we assume the Erd\"{o}s 
conjecture and deduce a much stronger consequence of it. We ask whether the sum
above can be arbitrarily large as the sets $A$ vary. Our first theorem answers
the question in the negative.

    Joseph L. Gerver \cite{Ger} proved that for every $\epsilon >0$, there exists
for all but a finite number of integers $ p \geq 3$, sets $S_p$ of positive integers, containing no arithmetic progression of $p$ terms, such that $\sum_{a\in S_p} \frac{1}{a} > (1-\epsilon)p ~\log p$. The set $S_p$ is the sequence $\{a_n \}$ where $a_1=1$ and for $n \geq 1$, $a_{n+1}$ is the smallest positive integer bigger than $a_n$ such that no $p$ elements of $a_1,a_2,\cdots a_{n+1}$in arithmetic progression.He guessed in that paper that for any prime $p$, the set $S_p$ may indeed
maximize the sum of the reciprocals of the elements of a set of positive integers having no p terms in arithmetic progression. On the other hand Joseph L. Gerver and L.Thomas Ramsey \cite{GR}
showed heuristically that the set $S_p$ is not maximizing the above sum for composite $p$. 
A corollary to our second theorem says that the Erd\" {o}s  conjecture implies the existence of a set of positive integers containing no $p$ terms in arithmetic progression which maximizes the above sum.  \\ [1mm]
In rest of the paper, $p$ is any fixed integer greater than or equal to 3.\\[4mm]

   The author thanks Professor Bhaskar Bagchi for valuable suggestions.
\section{Main Results}

\begin{thm}

  Let $\mathscr{A}_p$ be the collection of all subsets of $\mathbb{N}$ having
  no arithmetic progression of length $p$. Then, under the assumption of the
  Erd\"{o}s  conjecture, there is an absolute constant $B_p$ such that
  \begin{equation} \textnormal{Sup} \left\{ \sum_{a\in A} \frac{1}{a} : A \in \mathscr{A}_p \right\} \leq B_p . \end{equation}
\end{thm}

For further discussion, we need a topological structure on $\mathscr{A}_p$.
First we note that there is a natural one-to-one correspondence between the
power set $\mathscr{P}(\mathbb{N})$ and the set $\{ 0, 1 \}^{\mathbb{N}}$ of
all sequences of $0$s and $1$s ; namely, given any subset $A \subset
\mathbb{N}$, we send it to the sequence $\{ \delta_A ( n ) \}^{\infty}_{n =
1}$, where
$$ \delta_{A}(n) = \left\{\begin{array}{l}
     1 \textnormal{ if } n\in A\\
     0 \textnormal{ otherwise } .
   \end{array}\right. $$

Since $\{ 0, 1 \}^{\mathbb{N}}$ is compact by Tychonoff's theorem, the above
identification makes $\mathscr{P}(\mathbb{N})$ into a compact topological
space. In this topology, a sequence $\{ A_n \}$ of subsets converges to $A$
if, for any given $k$, there is some $N_k$ such that, whenever $n \geq
N_k$,
\begin{equation}
\label{conv}
  \delta_{A_n} ( j ) = \delta_A ( j ) \textnormal{ for } j = 1, 2, \cdots, k .
\end{equation}
Proposition 4 below says that $\mathscr{A}_p$ is a compact subspace of
$\mathscr{P}(\mathbb{N})$. For any set $A \in \mathscr{A}_p$, let us denote
the sum $\sum_{a\in A} \frac{1}{a}$ (which converges if we assume Erd\"{o}s conjecture) by
$\mu ( A )$. Then we have the following theorem.

\begin{thm}

  The map $A \longmapsto \mu ( A )$ between $\mathscr{A}_p$ and $[ 0, B_p ]$
  is continuous.
\end{thm}

Since $\mathscr{A}_p$ is compact, theorem 2 implies the following corollary.

\begin{cor}

  Under the assumption of the Erd\"{o}s  conjecture, there is a set $M_p \in
  \mathscr{A}_p$ such that
  \begin{equation} \mu ( X ) \leq \mu ( M_p ) \textnormal{ for all } X \in \mathscr{A}_p .
  \end{equation}
  That is, the supremum of the set $\{ \mu ( X ) : X \in \mathscr{A}_p \}$ is
  attained.
\end{cor}

\section{Proofs}

In this section, we shall present the proofs of theorem 1 and theorem 2. First we prove a proposition that will be needed later.

\begin{prop}
  $\mathscr{A}_p$ is a compact subspace of $\mathscr{P}(\mathbb{N})$.
\end{prop}

\begin{proof}

  Since $\mathscr{P}(\mathbb{N})$ is compact, it is enough to show that
  $\mathscr{A}_p$ is closed. Let $\{ A_n \}$ be sequence in $\mathscr{A}_p$
  converging to some $A \in \mathscr{P}(\mathbb{N})$. We need to show that $A
  \in \mathscr{A}_p$. Let us denote
  
  $A_n = \{ a_1^{( n )}, a_2^{( n )}, \cdots \}$ and $A = \{ a_1, a_2, \cdots
  \}$, where the terms in the sequences are written in the increasing order.
  Suppose, if possible, that $A \notin \mathscr{A}_p$. So there is an arithmetic
  progression $\{ a_{k_1}, a_{k_2}, \cdots, a_{k_p} \} \subset A$. We shall
  obtain a contradiction from this. Since $A_n \longmapsto A$, by the
  criterion (\ref{conv}) for convergence, we must have, for any given $k$, some integer
  $N_k$ such that,

\begin{equation} a_j^{( n )} = a_j \textnormal{ for } j = 1, 2, \cdots, k \end{equation}
  for all $n \geq N_k$. In particular, if $k = k_p$, we have, for $n \geq
  N_{k_p}$, 
  \begin{equation} a_{k_i}^{( n )} = a_{k_i} \textnormal{ for } i = 1, 2, \cdots, p. \end{equation}
  Since $\{ a_{k_i} : i = 1, 2, \cdots, p \}$ is an arithmetic progression,
  the above implies that $A_n \notin \mathscr{A}_p$ for $n \geq N_{k_p}$, which is a
  contradiction. So $A \in \mathscr{A}_p$ as was required to be proved.

\end{proof}

\vspace{6mm}

\noindent {\textbf{Proof of theorem 1}}

\begin{proof}

We shall prove this by contradiction. Let $A_0 = A \in \mathscr{A}_p$ be any finite set with $\sum_{a\in A} \frac{1}{a} = L> 0$.
 For example, we can take $A_0 = \{ 1 \}$. If we assume that the
statement of the theorem is not true, then we shall show that there is a
finite set $B \supset A$, $B \in \mathscr{A}_p$ with
\begin{equation} \sum_{b\in B} \frac{1}{b} \geq L + 1 . \end{equation}

 This will result in a contradiction to the conjecture of Erd\"{o}s  in the
following manner. Repeating this process that produces $B$ recursively, we get
an increasing sequence of sets $A_0 \subset A_1 \subset A_2 \subset \cdots$,
each of them finite and they all are in $\mathscr{A}_p$. Moreover, $$\sum_{a\in A_j}
\frac{1}{a} \geq L + j .$$

 Now the set $A_{\infty} = A_0 \cup A_1 \cup A_2 \cup
\cdots$ must be in $\mathscr{A}_p$ since any given collection of $p$ elements
in $A_{\infty}$ must also belong to $A_n$ for some $n$, so those elements can
not be in arithmetic progression. On the other hand, the sum $\sum_{a\in A}
\frac{1}{a}$ must diverge as it is bigger than any fixed number. So all that
is now left to prove the theorem is to produce such a set $B$, given $A$.

Let $N$ be the maximum of the elements of $A$. If the theorem is untrue, then
there must exist a set $E \in \mathscr{A}_p$ such that
\begin{equation} \sum_{e\in E} \frac{1}{e} \geq 2 N . \end{equation}
In fact, we may take $E$ to be a finite set; since, if $E$ is infinite, the
tail of the convergent sum will be small. Now define
\begin{equation} B = A \sqcup 2 {NE}, \end{equation}
where $\sqcup$ denotes disjoint union, and $2 {NE} = \{ 2 {Ne} : e
\in E \}$. Clearly $B$ is a finite set containing $A$, and
\begin{equation} \sum_{b\in B} \frac{1}{b} = \sum_{a\in A} \frac{1}{a} + \sum_{e\in E}\frac{1}{2 {Ne}} \geq L + 1 \end{equation}
by (7 , 8). Now to show that $B \in \mathscr{A}_p$, we first note that since $A \in
\mathscr{A}_p$ and $E \in \mathscr{A}_p$, no $p$ elements of either $A$ or $2
{NE}$ can be in arithmetic progression. Suppose, if possible, that $b_1,
b_2, \cdots, b_p \in B$ are in A.P., where $b_1, b_2, \cdots, b_k \in A$, and
$b_{k + 1}, b_{k + 2}, \cdots, b_p \in 2 {NE}$. If $k \geq 2$, then
\begin{equation} b_{k + 1} - b_k = b_k - b_{k - 1} . \end{equation}
Now, $b_k - b_{k - 1} \leq N - 1$ since $b_k, b_{k - 1} \in A$ and $N$ is the
maximum of the elements of $A$. But the right hand side, $b_{k + 1} - b_k \geq
b_{k + 1} - N \geq 2 N - N = N$, a contradiction. If $k = 1$, then
\begin{equation} b_2 - b_1 = b_3 - b_2, \end{equation}
or equivalently,
\begin{equation} b_1 = 2 b_2 - b_3 . \end{equation}
But $b_1 \leq N$, while $2 b_2 - b_3$ is a multiple of $2 N$ as both $b_2, b_3
\in 2 {NE}$. So we arrive at a contradiction again.  Hence we conclude that $B$ cannot have an arithmetic
progression of length $p$.

\end{proof}

For proving theorem 2, we first prove a lemma.

\begin{lem}

  Given any $\varepsilon > 0$, there exist a natural number $N$ such that for
  any $A \in \mathscr{A}_p$ with $\textnormal{ Min } A \geq N$,
  \begin{equation} \sum_{a\in A} \frac{1}{a} < \varepsilon . \end{equation}
\end{lem}

Note: In the above, $\textnormal{ Min } A$ denotes the smallest element in $A$.

\begin{proof}

Suppose, if possible, the lemma is not correct. Then there exists some $\varepsilon > 0$ such that for any given integer $M\geq 1$, there is a set $R \in \mathscr{A}_p$ depending on $M$ with the following properties:
\begin{equation} \label{mu} \mu ( R ) = \sum_{r\in R} \frac{1}{r} > \varepsilon, \end{equation}
and
\begin{equation} \label{r}
 \textnormal{Min } R \geq 2 M . \end{equation}
For that $\varepsilon$, we choose a set $A \in \mathscr{A}_p$ satisfying
\begin{equation}
\label{12}
 \mu ( A ) > M_p - \frac{\varepsilon}{12} . \end{equation}
where $M_p = \textnormal{Sup} \{ \mu ( A ) : A \in \mathscr{A}_p \} < \infty $ by Theorem 1. 
Let $A = \{ a_1, a_2, \cdots \}$ where $a_1 < a_2 < \cdots$ . Since $\sum_{a\in A}
\frac{1}{a} < \infty$, there is some $n_0$ such that
\begin{equation}
\label{12-2}
 \sum_{n = n_{0 + 1}}^{\infty} \frac{1}{a_n} < \frac{\varepsilon}{12} . \end{equation}
Let $A_1 = \{ a_1, a_2, \cdots, a_{n_0} \}$. Then
\begin{equation} 
\label{A1}
\mu ( A_1 ) > M_p - \frac{\varepsilon}{6} . \end{equation}
by (\ref{12}) and (\ref{12-2})

Now we take $M = a_{n_0}$ and write, $R = R_1 \sqcup R_2 \sqcup R_3 \sqcup R_4$ where
\begin{equation} R_j = R \bigcap \{ \bigsqcup_{i = 0}^{\infty} [ ( j + 1 ) 3^i M, ( j + 2 )
   3^i M ) \} ; j = 1, 2, 3, 4. \end{equation}
In other words, 
\begin{flushleft}
\begin{eqnarray}
 R_1 = R \bigcap \left\{ [ 2 M, 3 M ) \sqcup [ 6 M, 9 M ) \sqcup [ 18 M, 27
   M ) \sqcup \cdots \right\}, \nonumber\\ 
 R_2 = R \bigcap \left\{ [ 3 M, 4 M ) \sqcup [ 9 M, 12 M ) \sqcup [ 27 M, 36 M ) \sqcup \cdots \right\}, \nonumber\\
 R_3 = R \bigcap \left\{ [ 4 M, 5 M ) \sqcup [ 12 M, 15 M ) \sqcup [ 36 M,   45 M ) \sqcup \cdots \right\},\nonumber\\
 R_4 = R \bigcap \left\{ [ 5 M, 6 M ) \sqcup [ 15 M, 18 M ) \sqcup [ 45 M,  54 M ) \sqcup \cdots \right\}\nonumber . 
\end{eqnarray} 
\end{flushleft}
We have,
\begin{equation} \textnormal{Max } A_1 = M < 2 M \leq \textnormal{Min } R, \end{equation}
which implies $R \cap A_1 = \phi$, the empty set. Also, it is easy to check that no
$p$ elements of $A_1 \sqcup R_j, 1 \leq j \leq 4$, can be in an arithmetic
progression. So $A_1 \sqcup R_j \in \mathscr{A}_p$.

Since $\mu ( R ) > \varepsilon$, we must have
\begin{equation} \mu ( R_j ) > \frac{\varepsilon}{4} \end{equation}
for some $j, 1 \leq j \leq 4$.

For that $j$,
\begin{equation} \mu ( A_1 \sqcup R_j ) = \mu ( A_1 ) + \mu ( R_j ) > M_p +
   \frac{\varepsilon}{12}  \end{equation}
from (\ref{A1}). This is a contradiction to the fact that $M_p$ is the supremum of
the set $\{ \mu ( A ) : A \in \mathscr{A}_p \}$. This proves the lemma.

\end{proof}

\vspace{2mm}

Now we conclude this paper with the proof of theorem 2.

\vspace{2mm}

\noindent {\textbf{Proof of theorem 2}}

\begin{proof}
Suppose $\{ A_n \} \subset \mathscr{A}_p$ be a sequence and $A_n
\longrightarrow A$. We need to show that $\mu ( A_n ) \longrightarrow \mu ( A
)$.

Let us write the set $A$ as, $A = \{ a_1, a_2, a_3, \cdots \}$ where $a_1 <
a_2 < a_3 < \cdots$ and similarly for the sets $A_n$, we write
them as, $A_n = \{ a^{( n )}_1, a^{( n )}_2, \cdots \}$. Note
that if the set $A$ is finite, then $A_n = A$ for large enough $n$ and there
is nothing left to prove.  Let $\varepsilon > 0$ be any given real number. The lemma above
allows us to select an $N$ such that for any set $X \in \mathscr{A}_p$ with
$\textnormal{ Min } X \geq N$, we must have
\begin{equation}
 \label{epsi}
 \sum_{x\in X} \frac{1}{x} < \frac{\varepsilon}{2} .
 \end{equation}
Let $n_0$ be an integer such that $a_{n_0} \geq N$. Since $A_n \longrightarrow
A$, there is some $N_0$ such that $a^{( n )}_k = a_k$ for $1 \leq k \leq n_0$
and all $n \geq N_0$. Now, for $n \geq N_0$,

\begin{equation}
\begin{split}
 \Bigl\lvert \mu ( A_n ) - \mu ( A )\Bigr\rvert & = \Bigl\lvert  \sum_{k=n_0 + 1}^{\infty} \frac{1}{{a_k}^{(n)}} -
 \sum_{k=n_0 + 1}^{\infty} \frac{1}{a_k}  \Bigr \rvert \\
&\leq  \sum_{k=n_0 + 1}^{\infty} \frac{1}{{a_k}^{(n)}} + \sum_{k=n_0 + 1}^{\infty} \frac{1}{a_k}   < \varepsilon
\end{split}
\end{equation}
by (\ref{epsi}).
Hence $\mu(A_n) \longrightarrow \mu(A)$.

\end{proof}


\begin{thebibliography}{99}
\bibitem{Ger}Joseph L. Gerver, The sum of the reciprocals of a set of integers with no arithmetic progression of $k$ terms, Proceeding of American Mathematical Society, Vol. 62,
No. 2, Feb. 1977, pp. 211-214 .
\bibitem{GR}Joseph L. Gerver; L. Thomas Ramsey, Sets of Integers With No Long Arithmetic
Progression Generated by the Greedy Algorithm, Mathematics of Computation, Vol. 33, No.
148. (Oct., 1979.), pp.1353-1359.
\bibitem{GT}Ben Green; Terence Tao, arxiv.org/abs/math.NT/0404188
\bibitem{UL}Unpublished lecture, Faculte des Sciences, Paris, Dec. 4, 1975.
\end{thebibliography}
\end{document}